\documentclass[twocolumn,10pt,twoside,final,cleanfoot,cleanhead]{asme2ej}

\usepackage{epsfig} %% for loading postscript figures
\usepackage{graphicx}
\usepackage{subfigure}
\usepackage{mathtools}
\usepackage{algpseudocode}
\usepackage{algorithm}
\usepackage{amsmath}
\usepackage{bm}
\usepackage{amsfonts}
\usepackage{color}
\usepackage{verbatim}
\newcommand{\bA}{\mathbf{A}}
\newcommand{\bk}{\mathbf{k}}
\newcommand{\bx}{\mathbf{x}}

\newcommand{\by}{\mathbf{y}}
\newcommand{\bI}{\mathbf{I}}
\newcommand{\bK}{\mathbf{K}}

\newcommand{\btheta}{\boldsymbol{\theta}}
\newcommand{\bpsi}{\boldsymbol{\psi}}

\newcommand{\bxi}{\boldsymbol{\xi}}

\newcommand{\E}{\mathbb{E}}
\newcommand{\R}{\mathbb{R}}
\newcommand{\GP}{\operatorname{GP}}
\newcommand{\EEI}{\operatorname{EEI}}

\newcommand{\calN}{\mathcal{N}}

\newcommand{\sref}[1]{Sec.~\ref{sec:#1}}
\newcommand{\qref}[1]{Eq.~(\ref{eqn:#1})}

\newcommand{\fref}[1]{Fig.~\ref{fig:#1}}

\title{Extending Expected Improvement for High-dimensional
    Stochastic Optimization of Expensive Black-Box Functions}

\author{Piyush Pandita
        \affiliation{
          School of Mechanical Engineering\\
          Purdue University\\
          West Lafayette, Indiana 47907\\
        Email: ppandit@purdue.edu
        }   
        }

\author{Ilias Bilionis\thanks{Corresponding author} 
        \affiliation{School of Mechanical Engineering\\
        Purdue University\\
        West Lafayette, Indiana 47907\\
        Email: ibilion@purdue.edu
        }
        }

\author{Jitesh Panchal
        \affiliation{
        School of Mechanical Engineering\\
        Purdue University\\
        West Lafayette, Indiana 47907\\
        Email: panchal@purdue.edu
        }
        }

\begin{document}
\pagenumbering{arabic}
\maketitle    

\begin{abstract}{Design optimization under uncertainty is notoriously difficult 
when the objective function is expensive to evaluate. State-of-the-art techniques, 
e.g, stochastic optimization or sampling average approximation, fail to learn 
exploitable patterns from collected data and require an excessive number of 
objective function evaluations. There is a need for techniques that alleviate 
the high cost of information acquisition and select sequential simulations 
optimally. In the field of deterministic single-objective unconstrained global 
optimization, the Bayesian global optimization (BGO) approach has been relatively 
successful in addressing the information acquisition problem. BGO builds a 
probabilistic surrogate of the expensive objective function and uses it to 
define an information acquisition function (IAF) whose role is to quantify 
the merit of making new objective evaluations. Specifically, BGO iterates 
between making the observations with the largest expected IAF and rebuilding 
the probabilistic surrogate, until a convergence criterion is met. In this 
work, we extend the expected improvement (EI) IAF to the case of design 
optimization under uncertainty wherein the EI policy is reformulated to 
filter out parametric and measurement uncertainties.  To increase the 
robustness of our approach in the low sample regime, we employ a fully 
Bayesian interpretation of Gaussian processes by constructing a particle 
approximation of the posterior of its hyperparameters using adaptive 
Markov chain Monte Carlo. We verify and validate our approach by solving 
two synthetic optimization problems under uncertainty and demonstrate it 
by solving the oil-well-placement problem with uncertainties in the 
permeability field and the oil price time series.}
\end{abstract}

\section{Introduction}
\label{sec:intro}

The majority of stochastic optimization techniques are based on Monte Carlo
sampling, e.g., stochastic gradient descent \cite{bottou2010}, sample average
approximation \cite{kleywegt2002}, and random search \cite{shapiro2003}.
Unfortunately, the advantages offered by these techniques can be best leveraged
\cite{zinkevich2010} only when a large number of objective evaluations is
possible. Therefore, their applicability to engineering design/optimization
problems involving expensive physics-based models or even experimentally
measured objectives is severely limited.

Bayesian global optimization (BGO) has been successfully applied to the field
of single-objective unconstrained optimization.  \cite{torn1987, mockus1994,
locatelli1997, Jones2001, lizotte2008, benassi2011, bull2011}.  BGO builds a
probabilistic surrogate of the expensive objective function and uses it to
define an information acquisition function (IAF).  The role of the IAF is to
quantify the merit of making new objective evaluations.  Given an IAF, BGO
iterates between making the observation with the largest expected IAF and
rebuilding the probabilistic surrogate until a convergence criterion is met.
The most commonly used IAFs are the expected improvement (EI) \cite{Jones1998},
resulting in a version of BGO known as efficient global optimization (EGO), and
the probability of improvement (PoI) \cite{Jones2001}.  The operations research
literature has developed the concept of knowledge gradient (KG)
\cite{frazier2008,frazier2009,negoescu2011,scott2011}, which is essentially a
generalization of the EI, and the machine learning community has been
experimenting with the expected information gain (EIG)
\cite{villemonteix2009,hennig2012entropy,hernadez2014}.

BGO is not able to deal with stochastic optimization in a satisfactorily robust
way.  In this work, we propose a natural modification of the EI IAF, which is
able to filter out the effect of noise in the objective and, thus, enable
stochastic optimization strategies under an information acquisition budget.  We
will be referring to our version of EI as the Extended EI (EEI).  Our approach
does not suffer from the curse of dimensionality in the stochastic space, since
it represents both parametric and measurement noise in an equal footing and
does not explicitly try to learn the map between the uncertain parameters and
the objective.  However, we observed that naive applications of our strategy
fail to converge in the regime of low samples and high noise.  To deal with
this problem, we had to retain the full epistemic uncertainty of the underlying
objective surrogate.  This epistemic uncertainty corresponds to the fact that
the parameters of the surrogate cannot be determined exactly due to limited
data and/or increased noise.  Ignoring this uncertainty by picking specific
parameter values, e.g., by maximizing the marginal likelihood, typically yields
an overconfident, but wrong, surrogate.  This is a known problem in sequential
information acquisition literature, first mentioned by MacKay in
\cite{mackay1992}.  To avoid this issue, we had to explicitly characterize the
posterior distribution of the surrogate parameters by adaptive Markov chain
Monte Carlo sampling.  Remarkably, by keeping the full epistemic uncertainty
induced by the limited objective evaluations, we are able to characterize our
state of knowledge about the location of the optimum and the optimal value.

The outline of the paper is as follows. We start \sref{metho} by providing the
mathematical definition of the stochastic optimization problem that is being studied.
In \sref{gpr}, we introduce Gaussian process regression (GPR) which is used to
construct a probabilistic surrogate of the map between the design variables and
the objective.  In \sref{epistemic}, we show how the epistemic uncertainty on
the location of the optimum and the optimal value can be quantified.  In
\sref{EEI}, we derive our extension to EI suitable for stochastic optimization.
Our numerical results are presented in \sref{results}.  In particular, in
Sec.~\ref{sec:validate_1d} and~\ref{sec:validate_2d}, we validate our approach
using two synthetic stochastic optimization problems with known optimal
solutions and we experiment with various levels of Gaussian noise, as well as
heteroscedastic, i.e., input dependent, noise.  In \sref{owpp}, we apply our
methodology to solve the oil-well placement problem with uncertainties in soil
permeability and the oil price timeseries.  Our conclusions are presented in
\sref{conclusions}.

\section{Methodology}
\label{sec:metho}

We are interested in the following design optimization problem under uncertainty:
\begin{equation}
    \bx^{*} = \underset{\bx}{\arg\min}{\E_{\bxi}}\left[V(\bx;\bxi)\right],
    \label{eqn:optimal_design}
\end{equation}
where $V(\bx;\bxi)$ is the \emph{objective function} depending on a set of
\emph{design parameters} $\bx$ and \emph{stochastic parameters} $\bxi$.
The operator $\E_{\bxi}[\cdot]$ denotes the expectation over $\bxi$, i.e.,
\begin{equation}
    \E_{\bxi}\left[V(\bx;\bxi)\right] = \int V(\bx;\bxi)p(\bxi)d\bxi,
\end{equation}
where $p(\bxi)$ is the probability density function (PDF) of $\bxi$.
We will develop a methodology for the solution of \qref{optimal_design}
that addresses the following challenges:
\begin{enumerate}
    \item The objective is expensive to evaluate.
    \item It is not possible to compute the gradient of the objective with
          respect to $\bx$.
    \item The stochastic parameters $\bxi$ are either not observed directly, or
          they are 
          so high-dimensional that learning the dependence of the objective 
          with respect to them is impossible.
\end{enumerate}

Before we get to the specifics of our methodology, it is worth clarifying a
few things about the data collection process.
We assume that we can choose to evaluate the objective at any design point $\mathbf{x}$ we
wish. We envision this evaluation to take place as follows.
Behind the scenes, 
a random variable $\bxi$ is sampled from the, unknown, PDF $p(\bxi)$, and the
function $y=V(\bx;\bxi)$ is evaluated. We only see $y$ and not $\bxi$.
In this way, we can obtain an \emph{initial} data set consisting
of observed design points,
\begin{equation}
	{\bf{x}}_{1:n} = \{\it{\bf{x}}_{1},\cdots,\it{\bf{x}}_{n}\},
	\label{eqn:inputs}
\end{equation}
and the corresponding observed noisy objective evaluations,
\begin{equation}
	{\bf{y}}_{1:n} = \{\it{y}_{1},\cdots,\it{y}_{n}\}.
	\label{eqn:observed_outputs}
\end{equation}

What can be said about the solution of \qref{optimal_design}
using only the observed data $\bx_{1:n}$ and $\by_{1:n}$?
In the language of probability theory \cite{jaynes2003}, we would like to
characterize the probability of a design being optimal
conditional on the observations, and similarly for the optimal objective value.
Here probability corresponds to a state of belief and not to something random.
The uncertainty encoded in this probability is epistemic and it is 
induced by
the fact that inference is based on just $n$ observations. We will answer this
question by making no discounts on the Bayesian nature of Gaussian process
surrogates, see \sref{gpr} and \sref{epistemic}.

Where should we evaluate the objective next? Of course, looking for an optimal
    information acquisition policy is a futile task since the problem is
    mathematically equivalent to a non-linear stochastic dynamic programming
    problem~\cite{powell2012,bertsekas2007}. As in standard BGO, we will rely
    on a sub-optimal one-step-look-ahead strategy that makes use of an
    information acquisition function, albeit we will extend the EI information
    acquisition function so that it can cope robustly with noise, see \sref{EEI}.

\subsection{Gaussian process regression}
\label{sec:gpr}

Gaussian process regression \cite{rasmussen2006}
    is the Bayesian interpretation of classical Kriging 
\cite{cressie1990,Smith2014}. 
It is a powerful non-linear and non-parametric regression technique that 
has the added benefit of being able to quantify the epistemic uncertainties induced 
by limited data.
We will use it to learn the function that corresponds to
the expectation of the objective
$f(\cdot) = \E_{\bxi}[V(\cdot;\bxi)]$ from the observed data $\bx_{1:n}$ and
$\by_{1:n}$.

\subsubsection{Expressing prior beliefs}
\label{sec:prior}
A GP defines a probability measure on the space of meta-models, here
$f(\cdot)$, which can be used to encode our prior beliefs about the response,
e.g., lengthscales, regularity,
before we see any data.
Mathematically, we write:
\begin{equation}
    p(f( \cdot ) | \bpsi) = \GP(f(\cdot)|m(\cdot;\bpsi), k(\cdot,\cdot;\bpsi)),
\end{equation}
where $m(\cdot;\bpsi)$ and $k(\cdot,\cdot;\bpsi)$ are the mean and
covariance functions of the GP, respectively, and $\bpsi$ is a vector including
all the hyperparameters of the model.
Following the hierarchical Bayes framework, one would also have to specify
a prior on the hyperparameters, $p(\bpsi)$.

Note that information about the mean can actually be encoded in the covariance
function. Thus, without loss of generality, in this work we take $m(\cdot;\bpsi)$
to be identically equal to zero.
In our numerical examples, we will use the squared exponential (SE) covariance:
\begin{equation}
        k(\bx,\bx';\bpsi) = {s^2}\exp \left\{ { - \frac{1}{2}\sum\limits_{i = 1}^d 
	{\frac{{{{({x_i} - {x_i}')}^2}}}{{\ell_i^2}}} } \right\},
\end{equation}
where $d$ is the dimensionality of the design space, $s>0$ and $\ell_i>0$ can 
be interpreted as the signal strength of the response
and the lengthscale along input dimension $i$, respectively,
and $\bpsi  = \{ s,{\ell _1}, \ldots ,{\ell _d}\}$. Finishing, we assume that
all the hyperparameters are a priori independent:
\begin{equation}
    p(\bpsi) = p(s)\prod_{i=1}^dp(\ell_i),
    \label{eqn:psi_prior}
\end{equation}
where
\begin{equation}
    p(s) \propto \frac{1}{s}
    \label{eqn:s_prior}
\end{equation}
is the Jeffreys' prior \cite{jeffreys1946}, and
\begin{equation}
    p(\ell_i) \propto \frac{1}{1 + \ell_i^2}
\end{equation}
is a log-logistic prior \cite{conti2010}.

\subsubsection{Modeling the measurement process}
\label{sec:like}

To ensure analytical tractability, we assume that the measurement noise
is Gaussian with unknown variance $\sigma^2$.
Note that this could easily be relaxed to a student-t noise, which is more
robust to outliers.
The more general case of heteroscedastic, i.e., input-dependent, noise is an
open research problem and beyond the scope of the current work.
Note, however, that in our numerical examples we observe that our approach is
robust to modest heteroscedasticity levels.

Mathematically, the likelihood of the data is:
\begin{equation}
 p(\by_{1:n} | \bx_{1:n}, \btheta) = \calN\left(\by_{1:n}\middle|0, \bK_n(\bpsi) + \sigma^{2}\bI_n\right),
    \label{eqn:y_like}
\end{equation}
where $\calN(\cdot|\mu, \Sigma)$ is the PDF of a multivariate normal
random variable with mean $\mu$ and covariance matrix $\Sigma$,
$\bI_n\in\R^{n\times n}$ is the identity matrix,
$\bK_n(\bpsi) \in \R^{n\times n}$ is
the covariance matrix,
\begin{equation}
    \bK_n(\bpsi) = \left(\begin{array}{ccc}
        k(\bx_{1},\bx_{1};\bpsi) & \dots & k(\bx_{1},\bx_{n};\bpsi) \\
        \vdots	& 			 \ddots  &  \vdots\\
        k(\bx_{n},\bx_{1};\bpsi) & \dots & k(\bx_{n},\bx_{n};\bpsi)
    \end{array}\right),
\label{eq:covariance matrix}
\end{equation}
and, for notational convenience, we have defined $\btheta = \{\bpsi, \sigma\}$.
Finally, we need to assign a prior to $\sigma$.
We assume that $\sigma$ is a priori independent of all the variables in $\bpsi$
and set:
\begin{equation}
    p(\sigma) \propto \frac{1}{\sigma}.
    \label{eqn:sigma_prior}
\end{equation}

\subsubsection{Posterior state of knowledge}
\label{sec:posterior}
Bayes rule combines our prior beliefs with the likelihood of the data
and yields a posterior probability measure on the space of meta-models.
Conditioned on the hyperparameters $\btheta$, this measure is also a Gaussian process,
\begin{equation}
    \label{eqn:posterior_gp}
    p(f(\cdot)|\bx_{1:n}, \by_{1:n}, \btheta) = \GP\left(f(\cdot)\middle| m_n(\bx;\btheta), k_n(\bx,\bx';\btheta)\right),
\end{equation}
albeit with posterior mean and covariance functions,
\begin{equation}
    m_n(\bx;\btheta) = \left(\bk_n(\bx;\bpsi)\right)^{T}\left(\bK_n(\bpsi) + \sigma^2\bI_n\right)^{-1}\by_{1:n},
	\label{eqn:posterior_mean}
\end{equation}
and
\begin{eqnarray}
    k_n(\bx,\bx';\btheta)&=&k(\bx, \bx';\bpsi) \nonumber \\
                &&- \left(\bk_{n}(\bx;\bpsi)\right)^{T}\left(\bK_n(\bpsi) + \sigma^2\bI_N\right)^{-1} \bk_n(\bx';\bpsi) \nonumber\\
                &&
\label{eqn:predictive_covariance}
\end{eqnarray}
respectively,
where $\bk_n(\bx;\bpsi) = \left(k(\bx,\bx_1;\bpsi), \dots, k(\bx,\bx_n;\bpsi)\right)^T$,
and $\bA^T$  is the transpose of $\bA$.
Restricting our attention to a specific design point $\bx$,
we can derive from \qref{posterior_gp} the \emph{point-predictive probability
density} conditioned on the hyperparameters $\btheta$:
\begin{equation}
    \label{eqn:point_predictive}
    p(f(\bx)|\bx_{1:n}, \by_{1:n}, \btheta) = \calN\left(f(\bx)\middle|m_n(\bx;\btheta), \sigma_n^2(\bx;\btheta)\right),
\end{equation}
where $\sigma_n^2(\bx;\btheta) = k_n(\bx,\bx;\btheta)$.

To complete the characterization of the posterior state of knowledge, we need
to express our updated beliefs about the hyperparameters $\btheta$.
By a standard application of the Bayes rule, we get:
\begin{equation}
p(\btheta|\bx_{1:n},\by_{1:n}) \propto \
p(\by_{1:n}|\bx_{1:n},\btheta)p(\btheta),
\label{eqn:posterior_theta}
\end{equation}
where $p(\btheta) = p(\bpsi)p(\sigma)$. Unfortunately, \qref{posterior_theta}
cannot be computed analytically. Thus, we characterize it by a \emph{particle approximation}
consisting of $N$ samples, $\btheta_1,\dots,\btheta_N$ obtained by 
adaptive Markov chain Monte Carlo (MCMC)~\cite{haario2006}. Formally, we write:
\begin{equation}
    \label{eqn:theta_pa}
    p(\btheta|\bx_{1:n},\by_{1:n}) \approx \frac{1}{N}\sum_{i=1}^N\delta(\btheta - \btheta_i),
\end{equation}
where $\delta(\cdot)$ is Dirac's delta function. 
In our numerical results,
we use $N=90$ and the samples are generated as follows:
1) We obtain a starting point for the MCMC chain by maximizing the log of
the posterior \qref{posterior_theta}; 2) We burn  $10,000$ MCMC steps 
during which the MCMC proposal parameters are tuned; and 3) We perform
another
$90,000$ MCMC steps and record $\btheta$ every $1,000$ steps.

\subsection{Epistemic uncertainty on the solution of a stochastic optimization problem}
\label{sec:epistemic}
Now, we are in a position to quantify the epistemic uncertainty in the solution
of \qref{optimal_design} induced by the limited number of acquired data. 
Let $Q[\cdot]$ be any operator acting on functions $f(\cdot)$. Examples of 
such operators, are the minimum of $f(\cdot)$, $Q_{\min}[f(\cdot)] = \min_{\bx}f(\bx)$, or the
location of the minimum, $Q_{\arg\min}[f(\cdot)] = \arg\min_{\bx}f(\bx)$.
Conditioned on $\bx_{1:n}$ and $\by_{1:n}$ our state of knowledge about the value
of any operator $Q[\cdot]$ is
\begin{equation}
\begin{split}
    \label{eqn:Q_state}
    p(Q|\bx_{1:n},\by_{1:n}) = (\int (\int \delta\left(Q-Q[f(\cdot)]\right)p(f(\cdot)|\bx_{1:n}, \by_{1:n}, \btheta)\\
    df(\cdot))p(\btheta|\bx_{1:n},\by_{1:n})d\btheta),
\end{split}
\end{equation}
By sampling $M$ functions, $f_1(\cdot),\dots, f_M(\cdot)$ from 
\qref{posterior_gp} and using \qref{theta_pa}, we get the particle approximation:
\begin{equation}
    \label{eqn:Q_pa_state}
p(Q|\bx_{1:n},\by_{1:n}) \approx 
\frac{1}{NM} \sum_{i=1}^N\sum_{j=1}^M \delta\left(Q-Q[f_i(\cdot)]\right).
\end{equation}
 Our
derivation is straightforward and uses only the product and sum rules of probability
theory. The implementation, however, is rather technical.
For more details see the publications of Bilionis in the subject,
\cite{bilionis2012b,bilionis2012c,bilionis2013b,bilionis2014a,chen2015}.
In our numerical examples we use $M=100$.

\subsection{Extended expected improvement function}
\label{sec:EEI}

The classic definition of expected improvement, see \cite{Jones1998},
relies on the observed minimum $\tilde y_n = \min_{1\le i\le n} y_i$.
Unfortunately, this definition breaks down when $y_i$ is noisy.
To get a viable alternative, we have to filter out this noise.
To this end, let us define the
\emph{observed filtered minimum} conditioned on $\btheta$:
\begin{equation}
    \tilde m_n(\btheta) = \underset{1\le i\le n}\min m_n(\bx_i;\btheta),
\end{equation}
where $m_n(\bx;\btheta)$ is the posterior mean of \qref{posterior_mean}.
Using $\tilde m_n(\btheta)$, the improvement we would get if we observed 
$f(\bx)$ at design point $\bx$ is:
\begin{equation}
    I(\bx, f(\bx);\btheta) = \max\{0, \tilde m_n(\btheta) - f(\bx)\}.
\end{equation}
This is identical to the improvement function formulated in Sequential kriging 
optimization (SKO) \cite{Huang2006}. However, the EEI retains the full
epistemic uncertainty unlike SKO, which relies on a point estimate to the
hyper-parameters.
Since we don't know $f(\bx)$ or $\btheta$, we have to take their expectation
over our posterior state of knowledge, see \sref{posterior},
\begin{equation}
\begin{split}
    \EEI_n(\bx) = (\int \int I(\bx,f(\bx);\btheta)p(f(\bx)|\bx_{1:n},\by_{1:n},\btheta)df(\bx)\\
    p(\btheta|\bx_{1:n},\by_{1:n})d\btheta),
    \label{eqn:eei_def}
\end{split}
\end{equation}
where $p(f(\bx)|\bx_{1:n},\by_{1:n},\btheta)$ and $p(\btheta|\bx_{1:n},\by_{1:n})$
are given in \qref{point_predictive} and \qref{posterior_theta}, respectively.
The inner integral can be carried out analytically in exactly the same way as
one derives the classic expected improvement.
To evaluate the outer integral, we have to employ the particle approximation
to $p(\btheta|\bx_{1:n},\by_{1:n})$ given in \qref{theta_pa}.
The end result is:
\begin{equation}
\begin{split}
    \EEI_n(\bx) \approx \frac{1}{N}\sum_{i=1}^N[\sigma_n(\bx;\btheta_i)
            \phi\left (\frac{\tilde m_n(\btheta_i) - m_n(\bx;\btheta_i)}{\sigma_n(\bx;\btheta_i)}\right)\\
        + (\tilde m_n(\btheta_i) - m_n(\bx;\btheta_i))\Phi\left(\frac{\tilde m_n(\btheta_i) - m_n(\bx;\btheta_i)}{\sigma_n(\bx;\btheta_i)}\right)].
    \label{eqn:eei}
\end{split}
\end{equation}
Algorithm~\ref{alg:bgo} demonstrates how the derived information acquisition
criterion can be used in a modified version of BGO to obtain an approximation
to \qref{optimal_design}. Note that instead of attempting to maximize $\EEI_n(\bx)$
over $\bx$ exactly, we just search for the most informative point among a
set of $n_d$ randomly generated test points.
In our numerical examples we use $n_d=1,000$ test points following a
latin hypercube design~\cite{mckay2000}.

\begin{algorithm}[htb]
\caption{The Bayesian global optimization algorithm with the Extended expected improvement function}
\begin{algorithmic}[1]
\Require Observed inputs ${\bf{x}}_{1:n}$,
         observed outputs $\by_{1:n}$,
         number of candidate points tested for maximum EEI at each iteration $n_d$,
         maximum number of allowed iterations $S$,
         EEI tolerance  $\epsilon$. 
	\State $s \leftarrow 0$.
            \While {$s < S$} 
                \State Construct the particle approximation to the posterior of $\btheta$, \qref{theta_pa}.
                \State Generate a set of candidate test points $\hat\bx_{1:n_d}$, e.g., via a latin hypercube design \cite{mckay2000}.
                \State Compute EEI on all of the candidate points $\hat\bx_{1:n_d}$ using \qref{eei}.
                \State Find the candidate point $\hat\bx_j$ that exhibits the maximum EEI.
                \If{$\EEI_{n+s}(\bx_j) < \epsilon$}
                    \State Break.
                \EndIf
                \State Evaluate the objective at $\hat \bx_j$ measuring $\hat y$.
                \State $\bx_{1:n+s+1} \leftarrow \bx_{1:n+s}\cup\{\hat\bx_j\}$.
                \State $\by_{1:n+s+1} \leftarrow \by_{1:n+s}\cup\{\hat y\}$.
                \State $s\leftarrow s + 1$.
	\EndWhile
\end{algorithmic}
\label{alg:bgo}
\end{algorithm}

\section{Numerical Results}
\label{sec:results}

\begin{figure}[h!]
    \centering
    \subfigure[]{
        \includegraphics[width=0.5\textwidth]{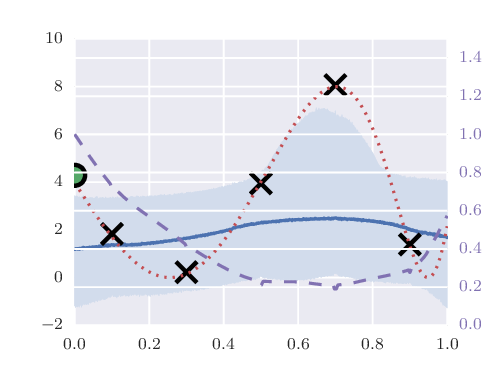}
    }
    \subfigure[]{
        \includegraphics[width=0.5\textwidth]{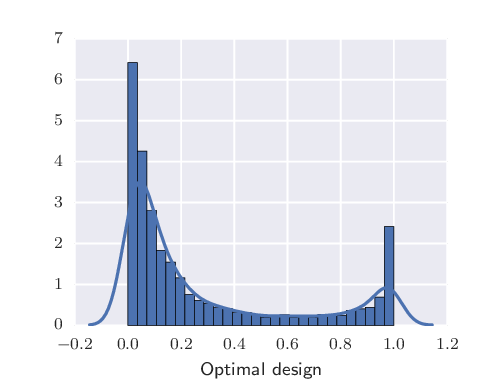}
    }
    \caption{One-dimensional synthetic example ($s(x)=0.1, n=5$).
        Subfigure~(a) depicts our initial state of knowledge about the true
        expected objective (dotted red line) conditioned on $n=5$ noisy
        observations (black crosses). 
        Subfigure~(b), shows a histogram of the predictive distribution of the
        optimal design $x^*$. 
    }
    \label{fig:ex1_1a}
\end{figure}

We validate our approach, see Sec.~\ref{sec:validate_1d} and~\ref{sec:validate_2d},
using two synthetic stochastic optimization problems with known optimal
solutions.
To assess the robustness of the methodology, we experiment with various levels
of Gaussian noise, as well as heteroscedastic, i.e., input dependent, noise.
In \sref{owpp}, we solve the oil-well placement problem with uncertainties in soil
permeability and the oil price timeseries.
Note that all the parameters required by our method, e.g., covariance function,
priors of hyperparameters, MCMC steps, have already been introduced in the
previous paragraphs and they are the same for all examples.
The only thing that we vary is the initial number of observations $n$.

\begin{figure}[h!]
    \centering
    \subfigure[]{
        \includegraphics[width=0.5\textwidth]{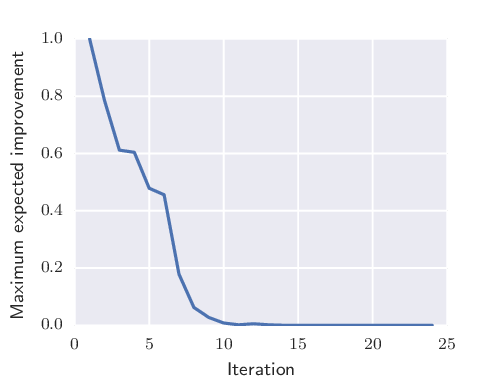}
    }
    \subfigure[]{
        \includegraphics[width=0.5\textwidth]{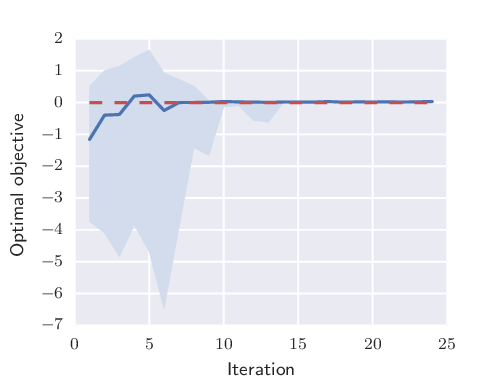}
    }
    \caption{One-dimensional synthetic example ($s(x)=0.1, n=5$).
        The dashed red line in Subfigure~(b) marks the real optimal value.
    }
    \label{fig:ex1_1b}
\end{figure}

\subsection{One-dimensional synthetic example}
\label{sec:validate_1d}
Consider the one-dimensional synthetic objective:
\begin{equation}
    V(x,\xi) = 4\left( {1 - \sin\left(6x + 8{e^{6x - 7}}\right)}\right) +
    s(x) \xi,
    \label{eqn:1d_example}
\end{equation}
for $x\in[0,1]$, where $\xi$ is a standard normal and for the noise standard deviation, $s(x)$,
we will experiment with $s(x)  = 0.01, 0.1, 1$, and
the heteroscedastic $s(x)  = \left(\frac{x-3}{3}\right)^{2}$.
Here, $\E_\xi[V(x,\xi)]$ is analytically available and it is quite trivial
to find that this function has two minima exhibiting the same objective value.

\begin{figure}[h!]
    \centering
    \subfigure[]{
        \includegraphics[width=0.5\textwidth]{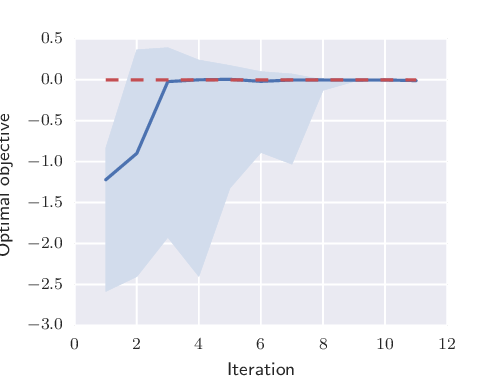}
    }
    \subfigure[]{
        \includegraphics[width=0.5\textwidth]{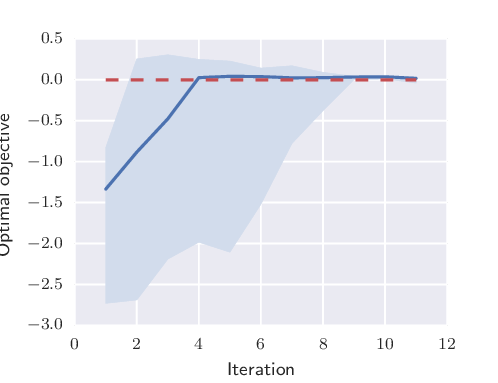}
   }
    \caption{One-dimensional synthetic example ($n=10$).}
    \label{fig:ex1_2a}
\end{figure}

\fref{ex1_1a}~(a) and~(b) visualize the posterior state of knowledge along with
the EEI (dashed purple line) as a function of $x$ and the epistemic uncertainty on the location of
the optimal design, respectively, for $s(x)=0.01$ when $n=5$. In \fref{ex1_1a}~(a), 
the solid blue line is the median of the predictive distribution of the GP and 
the shaded blue area corresponds to a $95\%$ prediction interval. 
\fref{ex1_1b}~(a) and~(b) depict the maximum EEI
and the evolution of the $95\%$ predictive bounds for the optimal objective
value (PBOO), respectively, as a function of the iteration number.
\fref{ex1_2a}~(a) and~(b) show the evolution of the PBOO for ($s(x)=0.01$) and ($s(x)=0.1$)
respectively and \fref{ex1_2a}~(a) and~(b) show the evolution of the PBOO 
for ($s(x)=1$) and ($s(x)=\left(\frac{x-3}{3}\right)^{2}$) respectively.

As expected, the larger the noise the more iterations are needed for
convergence. In general, we have observed that the method is robust to noise
as soon as the initial number of observations is not too low. For example,
the case $s(x)=1$ fails to converge to the truth, if one starts from less than
five initial observations.

\begin{figure}[h!]
    \centering
    \subfigure[]{
        \includegraphics[width=0.5\textwidth]{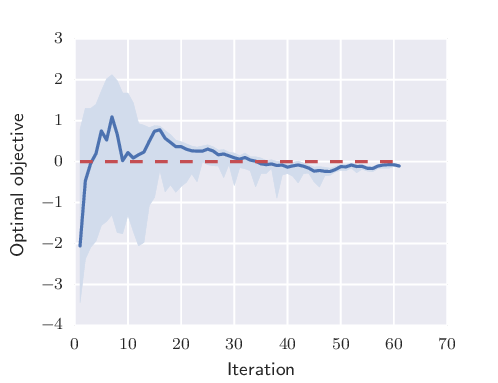}
    }
    \subfigure[]{
        \includegraphics[width=0.5\textwidth]{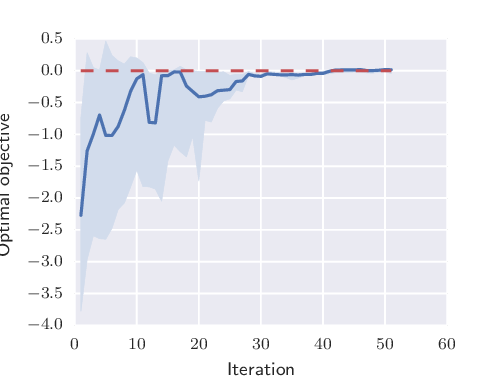}
    }
    \caption{One-dimensional synthetic example ($n=10$).}
    \label{fig:ex1_2b}
\end{figure}

% 1D example ends here, 2D toy example begins.
\subsection{Two-dimensional synthetic example}
\label{sec:validate_2d}
Consider the two-dimensional function \cite{sasena2002}:
\begin{equation}
\begin{split}
    V(\bx; \xi) = 2 + \frac{(x_{2}-x^{2}_{1})^{2}}{100} + (1-x_{1})^{2}  +2(2-x_{2})^{2}\\
+7\sin(0.5x_{2})\sin(0.7x_{1}x_{2}) + s(\bx) \xi,
\label{eqn:sasena}
\end{split}
\end{equation}
for $\bx\in[0,5]^2$, $\xi$ a standard normal,
and $s(\bx) = 0.01,0.1,1$, or the heteroscedastic $s(\bx) = (\frac{x_{2}-x{1}}{3})^{2}$.
As before, the expectation over $\xi$ is analytically available. It can easily
be verified that the objective exhibits three minima two of which are suboptimal.

\begin{figure}[htbp]
    \centering
    \subfigure[]{
        \includegraphics[width=0.50\textwidth]{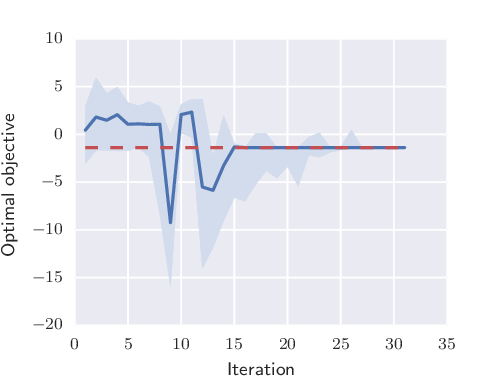}
    }
    \subfigure[]{
        \includegraphics[width=0.50\textwidth]{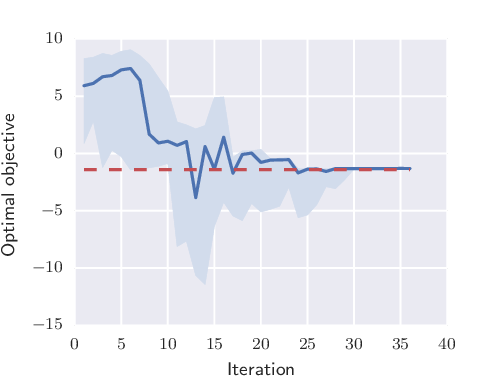}
    }
    \caption{Two-dimensional synthetic example ($n=20$).
        }
    \label{fig:ex2_1a}
\end{figure}

\fref{ex2_1a}~(a) and~(b) show the PBOO  for ($s(\bx)=0.01$) and ($s(\bx)=0.1$) 
and ~\fref{ex2_1b}~(a) and~(b) show the PBOO  for ($s(\bx)=1$) and 
($s(\bx)=\left(\frac{x_{2}-x{1}}{3}\right)^{2}$), respectively, as a function 
of the number of iterations. As before, the larger the noise the more 
iterations are required for convergence. 
The observed spikes are caused by the limited data used to build the surrogate.
In particular, the model is ``fooled" to believe that the noise is smaller than
it actually is and, as a result, it becomes more certain about the solution of
the optimization problem. As more observations are gathered though, the model
is self-corrected. This is a manifestation of the well known S-curve effect
of information acquisition~\cite[Ch. 5.2]{powell2012}. The existence of this
effect means, however, that one needs to be very careful in choosing the
stopping criterion.

\begin{figure}[htbp]
    \centering
    \subfigure[]{
        \includegraphics[width=0.50\textwidth]{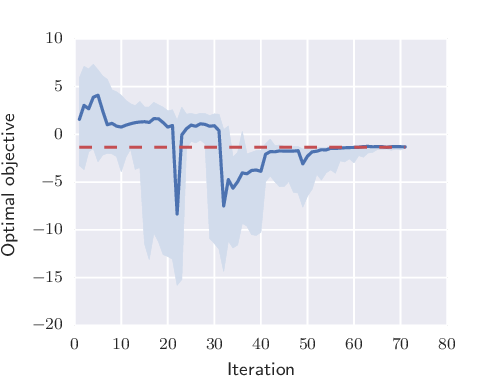}
    }
    \subfigure[]{
        \includegraphics[width=0.5\textwidth]{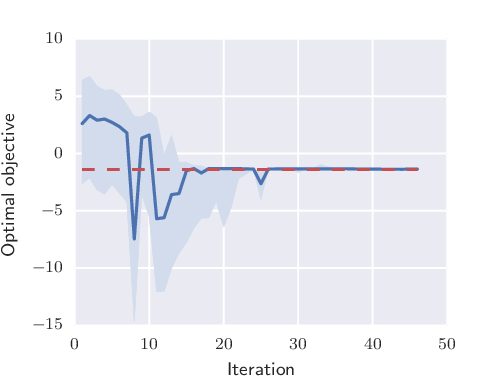}
    }
    \caption{Two-dimensional synthetic example ($n=20$).
        }
    \label{fig:ex2_1b}
\end{figure}
%Example 2D ends here. Oil well placement problem begins here.
\subsection{Oil well placement problem}
\label{sec:owpp}
During secondary oil production, water (potentially enhanced with chemicals or
gas) is injected into the reservoir through an \emph{injection} well. The
injected fluid pushes the oil out of the \emph{production} well. The \emph{oil}
\emph{well} \emph{placement} \emph{problem} (OWPP) involves the specification of the
number and location of the injection and production wells, the operating
pressures, the production schedule, etc., that maximize the net present value
(NPV) of the investment. This problem is of extreme importance for the oil
industry and an active area of research.  Several sources of uncertainty
influence the NPV, the most important of which are the time evolution of the
oil price (aleatoric uncertainty) and the uncertainty about the underground
geophysical parameters (epistemic uncertainty).

We consider an idealized 2D oil reservoir over the spatial domain
$\Omega = [0, 356.76]\times[0, 670.56]$ (measured in meters).
The four-dimensional design variable $\bx = (x_1, x_2, x_3, x_4)$
specifies the location of the injection well $(x_1,x_2)$, in which we pump water
(w), and the production well $(x_3,x_4)$,
out of which comes oil (o) and water.
Letting $\bx_s\in\Omega$ denote a spatial location, we assume that the permeability
of the ground is an isotropic tensor,
\begin{equation}
    \mathbf{C}(\bx_s;\bxi_c) = e^{g(\bx_s;\bxi_c)}c(\bx_s)\bI_3,
    \label{eqn:permeability}
\end{equation}
where $c(\bx_s)$ is the geometric mean (assumed to be the first layer
of the x-component of the SPE10 reservoir model permeability tensor~\cite{SPE10}),
$g(\bx_s;\bxi_c)$ is the truncated, at $13,200$ terms, Karhunen-Lo\`eve expansion of
a random field with exponential covariance function of lengthscale $\ell=10$ meters and
variance $10$, see \cite{ghanem2003}, and $\bxi_c$ is a ($13,200$)-dimensional
vector of standard normal random variables. Four samples of the permeability
field are depicted in \fref{permeability_samples}.
\begin{figure}[htbp]
    \centering
    \subfigure[]{
        \includegraphics[width=0.2\textwidth]{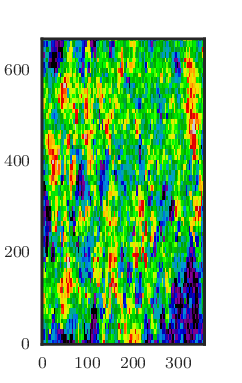}
    }
    \subfigure[]{
        \includegraphics[width=0.2\textwidth]{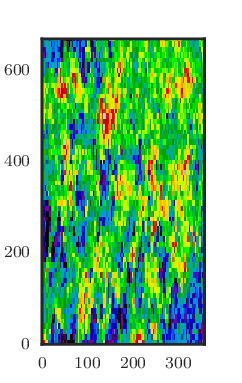}
    }
    \subfigure[]{
        \includegraphics[width=0.2\textwidth]{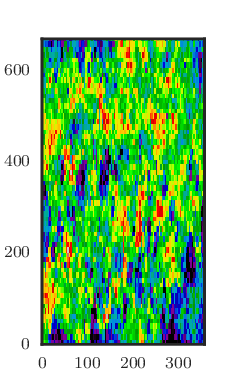}
    }
    \subfigure[]{
        \includegraphics[width=0.2\textwidth]{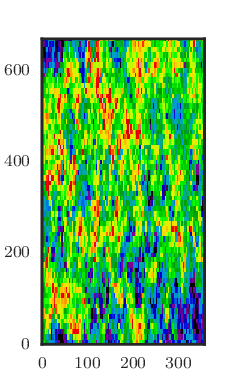}
    }
    \includegraphics[width=0.5\textwidth]{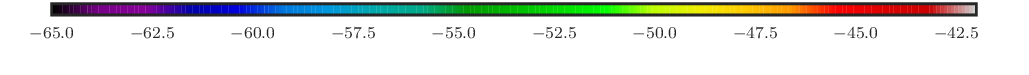}
    \caption{OWPP: Samples from the stochastic permeability model (in logarithmic scale)
    defined in \qref{permeability}.}
    \label{fig:permeability_samples}
\end{figure}

Given the well locations
$\bx$ and the stochastic variables $\bxi_c$, we solve a coupled system of time-dependent
partial differential equations (PDEs) describing the two-phase immiscible flow
of water and oil through the reservoir.
The solution is based on a finite volume scheme with a $60\times 220$ regular grid.
The form of the PDEs, the required boundary and initial conditions, as well as
the details of the finite volume discretization are discussed in~\cite{Bilionis2014}.
The parameters of the model that remain constant are as follows.
The water injection rate is $9.35\;\mbox{m}^3/\mbox{day}$, the connate water saturation
is $s_{\mbox{wc}} = 0.2$, the irreducible oil saturation is $s_{\mbox{or}} = 0.2$,
the water viscosity is set to $\mu_{w} = 3\times10^{-4}\;\mbox{Pa}\cdot\mbox{s}$, the oil viscosity to
$\mu_o = 3\times 10^{-3}\;\mbox{Pa}\cdot\mbox{s}$, the soil porosity is $10^{-3}$,
the timestep used is $\delta t = 0.1\;\mbox{days}$, and operations last
$T=2,000\;\mbox{days}$.
From the solution of the PDE system, we obtain the oil and water extraction rates
$q_o(t;\bx,\bxi_c)$ and $q_w(t;\bx,\bxi_c)$, respectively, where $t$ is the
time in days and the units of these quantities are in $\mbox{m}^3/\mbox{day}$.

\begin{figure}[htbp]
    \centering
    \includegraphics[]{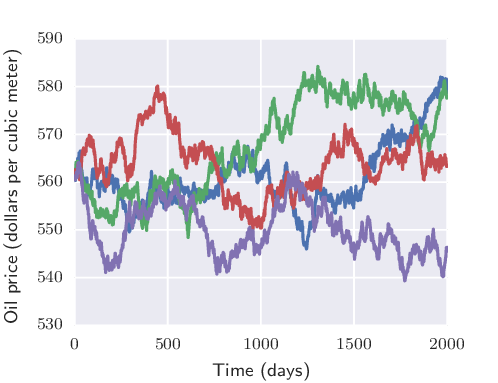}
    \caption{OWPP: Samples from the stochastic oil price model.}
    \label{fig:oil_price_samples}
\end{figure}

The oil price is modeled on a daily basis as $S_{o,t} = S_{o,0}e^{W_t}$,
where $S_{o,0} = \$560.8/\mbox{m}^3$, and $W_t$ is a random walk with a
drift:
\begin{equation}
    W_{t+1} = W_t + \mu  + \alpha\xi_{o,t},
\label{eqn:lognormal}
\end{equation}
where the $\mu = 10^{-8}$, $\alpha=10^{-3}$, and $\xi_{o,t}$ are independent
standard normal random variables. \fref{oil_price_samples} visualizes four
samples from the oil price model.
Since the process runs for
$T=2,000\;\mbox{days}$, we can think of $S_{o,t}$ as a function of the
$2,000$ independent identically distributed random variables $\bxi_o=\{\bxi_{o,1},\dots,\bxi_{o,T}\}$,
i.e., $S_{o,t} = S_{o,t}(\bxi_o)$.
For simplicity, we take the cost of disposing contaminated water is constant over
time $S_{w,t}^- = \$0.30/\mbox{m}^3$. Assuming a discount rate $r=10\%$ and risk neutrality, our objective is to
maximize the NPV of the investment. Equivalently, we wish to
minimize:
\begin{equation}
\begin{split}
    V(\bx;\bxi) = 10^{-6}\sum_{t=1}^{2,000\;\mbox{days}}\left[S_{w,t}q_(t;\bx,\bxi_c)-S_{o,t}(\bxi_o)q_o(t;\bx,\bxi_c)\right]\\
    (1+r)^{-t/365\;\mbox{days}},
\end{split}
\end{equation}
where $\bxi = \{\bxi_c, \bxi_o\}$, and the units are in million dollars.

\fref{ex3_results}~(a) shows the evolution of the PBOO as a function of the iterations of
our algorithm for the case of $n=20$ initial observations.
Note that in this case, we do not actually know what the optimal
value of the objective is. In subfigures~(b) and~(c) of the same figure, we 
visualize the initial set of observed well pairs and the well pairs selected
for simulation by our algorithm (where the blue `x' stands for the
injection well, the red `o' for the production well) respectively. Our algorithm quickly 
realizes the wells that are two close together are suboptimal and that it seems
to favor wells that are located at the bottom right and top right corners.
Note that the noise in this case is moderate, albeit heteroscedastic.

\begin{figure}[h!]
    \centering
    \subfigure[]{
        \includegraphics[width=0.5\textwidth]{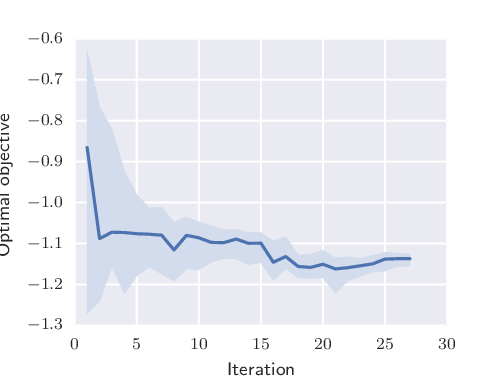}
    }
    \subfigure[]{
        \includegraphics[width=0.22\textwidth]{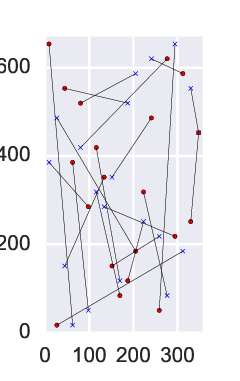}
    }
    \subfigure[]{
        \includegraphics[width=0.22\textwidth]{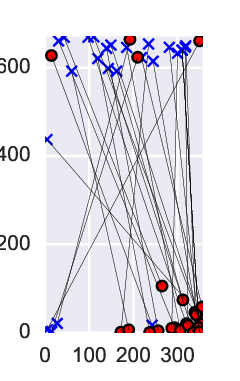}
    }
    \caption{OWPP ($n=20$).}
    \label{fig:ex3_results}
\end{figure}

\section{Conclusions}
\label{sec:conclusions}

We constructed an extension to the expected improvement which
makes possible the application of Bayesian global optimization to
stochastic optimization problems. In addition, we have shown how the epistemic
uncertainty induced by the limited number of simulations can be quantified, by
deriving predictive probability distributions for the location of the optimum
as well as the optimal value of the problem. We have validated our approach
with two synthetic examples with known solution and various noise levels, and
we applied it to the challenging oil well placement problem. 
The method offers a viable alternative to the sampling average approximation
when the cost of simulations is significant. We observe that our approach is
robust to moderate noise heteroscedasticity. There remain several open
research questions. In our opinion, the most important direction would be to
construct surrogates that explicitly model heteroscedasticity and use them
to extend the present methodology to robust stochastic optimization
and, subsequently, to multi-objective stochastic optimization.

\begin{acknowledgment}
Ilias Bilionis acknowledges the startup support provided by the School of 
Mechanical Engineering at Purdue University.
\end{acknowledgment}

\bibliographystyle{asmems4}

\bibliography{references} 

\end{document}